\documentstyle[11pt]{article}
\newtheorem{theorem}{Theorem}[section]
\newtheorem{lemma}[theorem]{Lemma}

\newtheorem{definition}[theorem]{Definition}

\newtheorem{claim}[theorem]{Claim}
\newtheorem{remark}[theorem]{Remark}
\begin{document}
\noindent{\centerline{\LARGE{
\bf On Extremal Elliptic $K3$ Surfaces\footnotemark[1]}} }
\input prepictex
\input pictex
\input postpictex
\vskip -0.5cm
\hskip -0.6cm   \hrulefill
\begin{flushright}

Qiang Ye 
 
{\em{ Departmant of Mathematics

National University of Singapore

 Kent Ridge 119260  

Singapore }}

{\it e-mail}: scip8188@nus.edu.sg

\end{flushright} 
\date{}

\footnotetext[1]{\noindent{Keywords}:\  Mordell-Weil Groups, \ Lattice.

\hskip 0.25cm \noindent{{\em Mathematics Subject Classification}}:
 14J28, 14J27.}
\section*{0 \ Introduction}
Let $C$ be a smooth projective curve over an algebraically
closed field of characteristic 0, and $X$ an {\em elliptic
surface} over $C$. By this we mean the following: $X$ is a 
smooth projective surface with a relatively minimal elliptic
fibration
$$
  f : X \longrightarrow C .
$$ 
In this paper, we also assume:

$(i)$\ $f$ has a global section $\cal O$, and

$(ii)$\ $f$ is not smooth, i.e., there is at 
least one singular fibre.
\vskip 0.5cm
\hskip -0.6cm To every (Jacobian) elliptic fibration $X$ there 
is a group of sections $\Phi (X)$ with the distinguished section 
$\cal O$  as zero. Up to a finite group, $\Phi (X)$ is identified 
with the relative automorphism group of the fibration.

\hskip -0.5cm Due to a formula of Shioda-Tate we have  the basic 
inequalities
$$
0\leq {\rm rank}\Phi \leq \rho (X)-2 
$$
where $\rho$ is the Picard number and the discrepancy
in the upper bound is related to the degree of reducibility
of the fibres.

\begin{definition}\hskip 0.5cm {\rm
An elliptic fibration $X$ is called {\em extremal} if and only
if $\rho (X)=h^{1,1} (X)$  ({\em maximal Picard number}) and
rank$\Phi (X)=0$.
}
\end{definition}
Let $f$ : $X \longrightarrow {\bf P}^1 $ be an {\em extremal} 
elliptic $K3$ surface. A fibration $f$ is called {\em semi-stable}
if each singular fiber of $f$ is of tpye $I_n$ [MP3]. 
Here we call a fibration {\em unsemi-stable} if it is not 
 semi-stable.
\vskip 0.5cm
\hskip -0.6cm In [MP2], R.Miranda and U.Persson have classified
possible semi-stable  fibrations.  The determination of all semi-stable 
fibrations has been done in [MP3] and [ATZ]. 
 In this paper, we first classify all possible configurations of
unsemi-stable fibrations (cf. Theorem 2.4). Then we calculate
the possible Mordell-Weil Groups for Case(A) (cf.Theorem 3.1),
 i.e., the case where each singular fibre of $f$ is {\em not}  
of tpye $I_n$. Finally, by using
the method in [ATZ], we will precisely 
determine which cases in {\em Table 1} are actually
realizable (cf.Theorem 0.4).
\vskip 0.5cm
\hskip -0.6cm Let $i_n$ denote the number of singular fibres of $f$ of type $I_n$. Similarly we define $i_n^*$, $ii$, $iii$, $iv$, $iv^*$, $iii^*$,
 $ii^*$ (cf. [MP1])  . Then we have the following Theorem 0.2.
\begin{theorem} Let 
$f$ : $X \longrightarrow {\bf P}^1 $ be an extremal elliptic
$K3$ surface with  $\deg J\not= 0$. Then
\[ 
\deg J = \sum_{n\geq 1} n(i_n + {i}^*_n )=
6\sum_{n\geq 1} (i_n + i^*_n ) + 4(ii + iv^* ) +
3(iii + iii^* ) + 2(iv + ii^* ) - 12.
\]
\end{theorem}
\begin{remark} {\rm  From [MP1,\ Lemma 3.1\ and 
\ Proposition 3.4], we know that the second equality in
{\rm Theorem\ 0.2} is replaced by ``$\leq$'' in general cases, 
for example, the fiber type $(II^*,II^*, I_2,I_1,I_1)$
{\rm (cf.[SI, Lemma 3.1])}.
 Thus it {\rm (Theorem\ 0.2)}justifies their naming ``extremal''. 
}
\end{remark}
\begin{theorem}
Let $f$ : $X \longrightarrow {\bf P}^1 $ be an extremal elliptic 
$K3$ surface and each singular fibre of $f$ is {\rm not} of
 type $I_n$. Then there exists exactly $11$ fiber types as
given below {\rm (table 1)}. In particular, Mordell-Weil Group is uniquely
determined by the fiber type of $f$. 
\vskip 0.5cm
\begin{center}
\begin{tabular}{| c|c|c|c|c|c|}\hline
 $\sharp$  &  the fibre type &  MW$(f)$ & $\sharp$  & the fibre type & MW$(f)$ 
\\ \hline

  $1$ & $(II^* ,I_1^* , I_1^* )$  &  $(0)$   &$ 7$ & $(IV^* ,IV^* , IV^* )$   &  ${\bf Z}/3{\bf Z}$       \\ \hline
 $ 2$ & $ (II^* ,II^* , IV)$ & $(0)$   &   $ 8$  &  $(IV^* ,IV^* , I_2^* )$     &  $(0)$   \\ \hline   
 $ 3$ & $(II^* ,IV^* , I_0^* )$ &  $(0)$     & $9$  &  $(IV^* ,I_3^* , I_1^* )$    &  $(0)$  \\ \hline 
$ 4$ & $(III^* ,III^* , I_0^* )$ &  ${\bf Z}/2{\bf Z}$     & $ 10 $  &  $(I_4^* ,I_1^* , I_1^* )$ & ${\bf Z}/2{\bf Z}$     \\ \hline 
 $ 5$ & $(III^* ,IV^* , I_1^* )$ & $(0)$  &  $11$ &$(I_2^* ,I_2^* , I_2^* )$ & ${\bf Z}/2{\bf Z}\oplus {\bf Z}/2{\bf Z} $     \\ \hline 

$ 6$ & $(III^* ,I_2^* , I_1^* )$ &  ${\bf Z}/2{\bf Z}$  &  & & \\ \hline  
\end{tabular} 
\end{center} 
All the above $11$ fiber types are realizable.
\end{theorem}
This paper is organized as follows.
In Section 1, we introduce some basic notation and theorems
 which will be used
in the paper. In Section 2, we first prove Theorem 0.2.
 Then we  give the combinatorical
classification of the possible unsemi-stable fibration 
(cf.Theorem 2.4). In Section 3, we calculate  all possible Mordell-Weil
Groups for Case(A) (cf. Theorem 3.1).
In Section 4, we prove Theorem 0.4.

\hskip -0.6cm At the same time, and indepedently, 
I.Shimada and D.Q.Zhang  present
a complete list of extremal elliptic $K3$ surfaces by using the
different method [SZ].
\vskip 0.5cm
\hskip -0.6cm {\bf Acknowledgement.} I would like to thank my
advisor Professor D.-Q. Zhang for introducing me to this
subject, lending his precious manuscript to me and
many enlightening instructions during my preparation of this paper.
 I also wish to express my sincere gratitude to
Professor J.Conway and Professor N.Sloane for 
their helpful correspondences.
\section{Preliminaries}

{\bf (a)  Lattices}
\vskip 0.5cm
\hskip -0.5cm Let $L$ be a lattice, i.e.,

$(i)$ \ $L$ is a free finite {\bf Z} module and

$(ii)$ \ $L$ is equipped with a non-degenerate bilinear
symmetric pairing $<  ,  >$.

The determinant of $L$, det$L$, is defined as
 the determinant of the matrix $I=(<x_i , x_j >)$
where $\{ x_1 ,\dots ,x_r \}$ is a {\bf Z}-basis of $L$
($r$= the rank of $L$):
$$
  {\rm det}L = {\rm det}(<x_i ,x_j >).
$$
We define the positive- (or negative-) definiteness or
the signature of a lattice by that of the matrix $I$,
noting that these properties
 are independent of the choice of a basis. An lattice $L$ is
 called {\em even} if $<x,x> \in 2{\bf Z}$ for all $x \in L$.
 We call $L$ {\em unimodular} if det$L$=1. Let $J$ be a
 sublettice of $L$. We denote its orthogonal complement with 
respect to $< , >$ by $J^{\perp}$.
 
For a lattice $L$, we denote its {\em dual} lattice by $L^{\vee}$. By
using pairing, $L$ is embedded in $L^{\vee}$ as a sublettice with 
the same rank. Hence the quotient group $L^{\vee}/L$ is a finite
abelian group, which we denote by $G_J$.

For an even lattice $L$, we define a quadratic form $q_L$ with
values in {\bf Q}/2{\bf Z} as follows:
$$
q_L (x \bmod{L}) = <x, x> \bmod{ 2{\bf Z}}.
$$
\begin{lemma}  For $j=1$,$2$, let 
$\Delta_j = \Delta (1)_j \oplus \cdots \oplus \Delta (r_j )_j$
be a lattice where each $\Delta (i)_j $is of Dynkin type $A_a$,
$D_d$ or $E_e$.

$(1)$ \ Suppose that $\Phi$ : $\Delta_1 \longrightarrow \Delta_2$
is a lattice-isometry. Then $r_1 =r_2$ and $\Phi (\Delta (i)_1 )=
\Delta (i)_2$ after relabelling.

$(2)$ \ Let ${\bf B}(6)= E_7 \oplus D_6 \oplus D_5$,  
${\bf B}(10)= D_8 \oplus D_5 \oplus D_5$,
${\bf B}(11)= D_6 \oplus D_6 \oplus D_6$.
Then we have

$(i)$\ ${\bf B}(6) \subset E_7 \oplus D_{11}$ is an  index-$2$ 
lattice extension.

$(ii)$\ ${\bf B}(10) \subset D_5 \oplus D_{13}$ is an index-$2$ 
extension, ${\bf B}(10) \subset D_8 \oplus D_{10}$ is 
an index-$2$ extension and  ${\bf B}(10) \subset D_{18}$ is 
an index-$4$ extension.

$(iii)$\ ${\bf B}(11) \subset D_6 \oplus D_{12}$ is an index-$2$ 
extension, ${\bf B}(11) \subset D_{18}$ is an index-$4$ extension.
\end{lemma}
{\em Proof.}\hskip 0.5cm 
We observe that
$$
 |\det (A_n )|= n+1 , |\det (D_n )|=4 , |\det (E_6 )|=3 ,
|\det (E_7 )|=2 , |\det (E_8 )|=1,
$$
and for an index-$n$ lattice extension $L \subset M$ one has
$$
  |\det (L)| = n^2 |\det (M)|.
$$
Then (1) comes  from [ATZ, Lemma 1.3], 
and  (2) can be obtained by an easy calculation.
\begin{definition} {\rm (The lattice $D_n$) [CS]
 \  For $n\geq 3$,
$$
D_n = \{(x_1,x_2,...,x_n)\in {\bf Z}^n : x_1+\cdots +x_n \ is\ even \}.
$$
}
\end{definition}
\begin{remark} {\rm From {\rm Lemma 1.1}, we know 
$D_8 \oplus D_5 \oplus D_5 \subset D_{18}$ is an index-$4$
extension. Thus $D_{18}/(D_8 \oplus D_5 \oplus D_5)$ maybe
${\bf Z}/4{\bf Z}$ or
${\bf Z}/2{\bf Z} \oplus  {\bf Z}/2{\bf Z}$. In the
following {\rm Lemma 1.4}, We shall prove that,
 $D_{18}/(D_8 \oplus D_5 \oplus D_5)
={\bf Z}/2{\bf Z} \oplus  {\bf Z}/2{\bf Z} $.  
}
\end{remark}

\begin{lemma} {\it For any lattice-isometric embedding 
$i$: $ D_5 \oplus D_5 \oplus D_8 \longrightarrow D_{18}$, we have
$$
D_{18} / (D_5 \oplus D_5 \oplus
D_8) = {\bf Z}/2{\bf Z} \oplus  {\bf Z}/2{\bf Z}. 
$$
}
\end{lemma}
{\em Proof.}\hskip 0.5cm 
We denote $A\oplus B\oplus C = D_5 \oplus D_5 \oplus D_8$ and let
$i$ : $A \oplus B \oplus C \longrightarrow  D_{18}$ 
be a lattice-isometric embedding. For one generator
$e$ of the lattices $A$, $B$ or $C$, we  assume that 
$i(e)=(x_1 ,x_2 ,...,x_{18}) \in D_{18}$. 
Since $2=<e,e>=<i(e),i(e)>=\sum_{i=1}^{18} x_i^2$ and
$x_i$ is integer, we have

\begin{claim}\hskip 0.5cm  {\it There are exactly two
 coordinates of $i(e)$  which are non-zero and each of which is $1$ or $-1$. 
}
\end{claim}
Thus by relabelling the coordinates, we may assume that one
 generators  $e_1$ of $A$ satisfies
$$
 i(e_1 )=(1,1,0,...,0) \ or \ (1,-1,0,...,0).
$$
Then we can use the connections among the generators
 in Dynkin diagram of $D_5$ to get the possible coordinates of
the generators $e_1$, $e_2$, $e_3$, $e_4$ , $e_5$ of $A$ .
After a simple calculation, we find that, by 
rebelling the coordinates, we may assume that 
$$
i(A) \subset (x_1,x_2,x_3,x_4,x_5,0,...,0) \cap D_{18} := L_1.
$$
On the other hand, we know $L_1$ is a $D_5$ type lattice and
$i$ is  lattice isometry, thus we get
$$
i(A) = L_1.
$$ 
By using the same method,  we may assume that
$$
i(B)=(0,..,0,x_6,x_7,..,x_{10},0,..,0) \cap D_{18} := L_2
$$
and
$$
i(C)=(0,..,0,x_{11},,..,x_{18}) \cap D_{18} := L_3.
$$
Here we will use the orthonormal conditions among $i(A)$,
$i(B)$ and $i(C)$.

A direct computation shows that
$$
D_{18}/ (L_1 \oplus L_2 \oplus L_3) = 
{\bf Z}/2{\bf Z} \oplus {\bf Z}/2{\bf Z}.
$$ 
Thus we prove the Lemma 1.4.
\vskip 0.5cm
\hskip -0.6cm By using the same idea as above proof, we get

\begin{theorem}
 {\it For $m=\sum_{i=1}^k n_i$, we have
$$
D_{m} /  (\oplus_{i=1}^k  D_{n_i} ) = \oplus_{i=1}^{k}  {\bf Z}/2{\bf Z} 
$$
for any lattice-isometric embedding 
$i$ $:$ $\oplus_{i=1}^k  D_{n_i} \longrightarrow
D_m$.
}
\end{theorem}
{\bf (b) \ Mordell-Weil lattices of elliptic surface} 
\vskip 0.5cm
\hskip -0.6cm  Given  an ellipric surface $f$: $X\longrightarrow C$,
let $F_{\nu} = f^{-1} ({\nu})$ denote the fibre over $\nu \in
C$, and let

$  Sing(f) =\{ \nu \in C | F_{\nu} \ {\rm is\  singular} \}$.

$  {\bf R}= {\rm Red}(f) = \{ \nu \in C | F_{\nu}\  {\rm is \ 
reducible}\}. $

\hskip -0.6cm For each $\nu \in {\bf R}$, let
$$
  F_{\nu} = f^{-1} (\nu ) = \Theta_{\nu ,0} + \sum_{i=1}^{m_{\nu}
-1} \mu_{\nu ,i} \Theta_{\nu ,i}  \ \ (\mu_{\nu ,0} =1 )
$$
where $ \Theta_{\nu ,i}$ ($ 0\leq i \leq m_{\nu} -1$) are the 
irreducible components of $F_{\nu}$, $m_{\nu}$ being their
number, such that $ \Theta_{\nu ,0}$ is the unique component
of $F_{\nu}$ meeting the zero section.

\hskip -0.6cm Here we denote
$$
  E(K) = {\rm the \ group\ of \ sections \ of\ }\  f,
$$
and
$$
  NS(X) = {\rm the \ group \ of \ divisors \ on}\  X \ \ {\rm
 modulo \ algebraic \ equivalence.}
$$
\begin{theorem} {\rm (cf.[Sh, Theorem 1.1,1.2,1.3])}
Under the assumptions for the elliptic surfaces in Introduction, 
we have 

$(1)$\ $E(K)$ is a finite generated abelian group.

$(2)$\ $N(X)$ is finitely generated and torsion-free.

$(3)$\ Let $T$ denote the subgroup of $NS(X)$ generated by
the zero section $({\cal O})$ and all the irreducible 
components of fibres. Then, there is a natural isomorphism
$$
  E(K) \cong NS(X)/T ,
$$
which maps $P\in E(K)$ to $(P)$ mod $T$.
\end{theorem}
\begin{theorem} {\rm (cf.[Sh, Lemma 8.1])}
For any $P$, $Q \in E(K)$, let
$$
  <P,Q> = - (\varphi (P) \cdot \varphi (Q))   \hskip 2cm (*)
$$
where $\varphi (P)$ (resp.  $\varphi (Q)$), satisfying
the condition:

$(1)$ \ $\varphi (P) \equiv (P)$ {\rm mod} $T_{\bf Q}$, and

$(2)$ \  $\varphi (P) \perp T$.

Then it defines a symmetric bilinear pairing on $E(K)$, which
induces the structure of a positive-definite lattice on
$E(K)/E(K)_{\it tor}$.
\end{theorem}

\begin{definition} {\rm 
The pairing $(*)$ on the Mordell-Weil group $E(K)$ is called
the {\em height pairing}, and the lattice
$$
   (E(K)/E(K)_{\it tor} , <,>)
$$
is called the {\em Mordell-Weil Lattice} of the elliptic
curve $E/K$ or of the elliptic surface $f$: $S \longrightarrow
C$.
}
\end{definition}

\begin{theorem} {\rm (Explicit formula for the height pairing) 
[Sh, Theorem 8.6]} For any $P$,$Q \in E(K)$, 
we have
$$
<P ,Q> = \chi + (P{\cal O}) + (Q{\cal O}) - (PQ) - \sum_{\nu \in R} contr_{\nu}
(P,Q),
$$
$$
<P,P>=2\chi + 2(P{\cal O}) - \sum_{\nu \in R} contr_{\nu}
(P).
$$
\end{theorem}
\begin{remark} {\rm Here $\chi$ is the arithmetic genus of 
$S$, and $(P{\cal O})$ is the intersection number of the sections $(P)$
and $({\cal O})$, and similarly for $(Q{\cal O})$,$(PQ)$. The term  $contr_{\nu}
(P,Q)$ stands for the local contribution ar $\nu \in R$, which is
 defined as follows: suppose that $(P)$ interests $\Theta_{\nu ,i}$
and $(Q)$ intersects $\Theta_{\nu ,j}$. Then we let
$$
contr_{\nu} (P,Q) = \left \{\begin{array}{ll} 
                (-A_{\nu}^{-1})_{i,j} ,& if \ i\geq 1,j\geq 1, \\
                               0                 ,& otherwise.
 \end{array}
\right.
$$
where the first one means the $(i,j)$-entry of the matrix $(-A_{\nu}^{-1})$.
Further we set
$$
 contr_{\nu} (P) =  contr_{\nu} (P,P).
$$

\hskip -0.5cm Arrange $\Theta_i = \Theta_{\nu ,i}$ $(i=0,1,\cdots ,m_{\nu} -1)$ so 
that the simple components are numbered as in the figure below.
\vskip 0.5cm
\centerline{\font\thinlinefont=cmr5
\begingroup\makeatletter\ifx\SetFigFont\undefined
% extract first six characters in \fmtname
\def\x#1#2#3#4#5#6#7\relax{\def\x{#1#2#3#4#5#6}}%
\expandafter\x\fmtname xxxxxx\relax \def\y{splain}%
\ifx\x\y   % LaTeX or SliTeX?
\gdef\SetFigFont#1#2#3{%
  \ifnum #1<17\tiny\else \ifnum #1<20\small\else
  \ifnum #1<24\normalsize\else \ifnum #1<29\large\else
  \ifnum #1<34\Large\else \ifnum #1<41\LARGE\else
     \huge\fi\fi\fi\fi\fi\fi
  \csname #3\endcsname}%
\else
\gdef\SetFigFont#1#2#3{\begingroup
  \count@#1\relax \ifnum 25<\count@\count@25\fi
  \def\x{\endgroup\@setsize\SetFigFont{#2pt}}%
  \expandafter\x
    \csname \romannumeral\the\count@ pt\expandafter\endcsname
    \csname @\romannumeral\the\count@ pt\endcsname
  \csname #3\endcsname}%
\fi
\fi\endgroup
\mbox
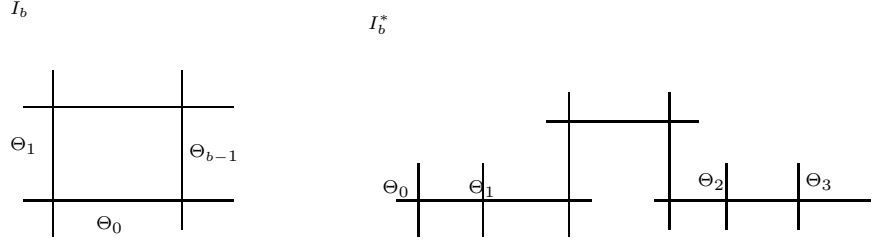
{\beginpicture
\setcoordinatesystem units <0.60000cm,0.60000cm>
\unitlength=0.60000cm
\linethickness=1pt
\setplotsymbol ({\makebox(0,0)[l]{\tencirc\symbol{'160}}})
\setshadesymbol ({\thinlinefont .})
\setlinear
%
% Fig POLYLINE object
%
\linethickness= 0.500pt
\setplotsymbol ({\thinlinefont .})
\putrectangle corners at  1.746 18.415 and  6.350 18.415
%
% Fig POLYLINE object
%
\linethickness= 0.500pt
\setplotsymbol ({\thinlinefont .})
\putrectangle corners at  1.905 16.351 and  6.350 16.351
%
% Fig POLYLINE object
%
\linethickness= 0.500pt
\setplotsymbol ({\thinlinefont .})
\putrectangle corners at  1.746 16.351 and  1.905 16.351
%
% Fig POLYLINE object
%
\linethickness= 0.500pt
\setplotsymbol ({\thinlinefont .})
\putrectangle corners at  2.381 19.209 and  2.381 15.558
%
% Fig POLYLINE object
%
\linethickness= 0.500pt
\setplotsymbol ({\thinlinefont .})
\putrectangle corners at  5.239 19.209 and  5.239 15.716
%
% Fig POLYLINE object
%
\linethickness= 0.500pt
\setplotsymbol ({\thinlinefont .})
\putrectangle corners at 10.001 16.351 and 14.287 16.351
%
% Fig POLYLINE object
%
\linethickness= 0.500pt
\setplotsymbol ({\thinlinefont .})
\putrectangle corners at 13.335 18.098 and 16.669 18.098
%
% Fig POLYLINE object
%
\linethickness= 0.500pt
\setplotsymbol ({\thinlinefont .})
\putrectangle corners at 15.716 16.351 and 20.479 16.351
%
% Fig POLYLINE object
%
\linethickness= 0.500pt
\setplotsymbol ({\thinlinefont .})
\putrectangle corners at 13.811 18.733 and 13.811 15.558
%
% Fig POLYLINE object
%
\linethickness= 0.500pt
\setplotsymbol ({\thinlinefont .})
\putrectangle corners at 16.034 18.733 and 16.034 15.716
%
% Fig POLYLINE object
%
\linethickness= 0.500pt
\setplotsymbol ({\thinlinefont .})
\putrectangle corners at 10.478 17.145 and 10.478 15.558
%
% Fig POLYLINE object
%
\linethickness= 0.500pt
\setplotsymbol ({\thinlinefont .})
\putrectangle corners at 11.906 17.145 and 11.906 15.558
%
% Fig POLYLINE object
%
\linethickness= 0.500pt
\setplotsymbol ({\thinlinefont .})
\putrectangle corners at 17.304 17.145 and 17.304 15.716
%
% Fig POLYLINE object
%
\linethickness= 0.500pt
\setplotsymbol ({\thinlinefont .})
\putrectangle corners at 18.891 17.145 and 18.891 15.716
%
% Fig POLYLINE object
%
\linethickness= 0.500pt
\setplotsymbol ({\thinlinefont .})
\putrectangle corners at 10.478 17.145 and 10.478 16.828
%
% Fig POLYLINE object
%
\linethickness= 0.500pt
\setplotsymbol ({\thinlinefont .})
\putrectangle corners at 10.478 16.034 and 10.478 15.558
%
% Fig TEXT object
%
\put{\SetFigFont{7}{8.4}{rm}$I_b$} [lB] at  1.429 20.479
%
% Fig TEXT object
%
\put{\SetFigFont{7}{8.4}{rm}  $\Theta_1$ } [lB] at  1.270 17.462
%
% Fig TEXT object
%
\put{\SetFigFont{7}{8.4}{rm}$\Theta_0$} [lB] at  3.334 15.716
%
% Fig TEXT object
%
\put{\SetFigFont{7}{8.4}{rm}$\Theta_{b-1}$} [lB] at  5.397 17.304
%
% Fig TEXT object
%
\put{\SetFigFont{7}{8.4}{rm}$\Theta_0$} [lB] at  9.684 16.510
%
% Fig TEXT object
%
\put{\SetFigFont{7}{8.4}{rm}$\Theta_1$} [lB] at 11.589 16.510
%
% Fig TEXT object
%
\put{\SetFigFont{7}{8.4}{rm}$\Theta_2$} [lB] at 16.669 16.669
%
% Fig TEXT object
%
\put{\SetFigFont{7}{8.4}{rm}$\Theta_3$} [lB] at 19.050 16.669
%
% Fig TEXT object
%
\put{\SetFigFont{7}{8.4}{rm}$I_b^*$} [lB] at  9.366 20.161
\linethickness=0pt
\putrectangle corners at  1.270 20.892 and 20.504 15.532
\endpicture}
}
\vskip 0.5cm
\hskip -0.6cm For the other types of reducible fibres, the numbering is irrelevant.
Assume that $(P)$ intersects $\Theta_{\nu ,i}$ and $(Q)$ intersect
$\Theta_{\nu ,j}$ with $i>1$,$j>1$. Then we have the following 
table: the forth row is for the case $i<j$ (interchange $P$, $Q$
if necessary).
\vskip 0.5cm
\begin{center}
\begin{tabular}{|c|c|c|c|c|c|c|}\hline
$T_{\nu}^-$ & $A_1$ & $E_7$ & $A_2$ & $E_6$ & $A_{b-1}$ &  $D_{b+4}$ 
\\ \hline 
type of $F_{\nu}$ & $III$ & $III^*$ & $IV$ & $IV^*$ & $I_b (b\geq 2 )$
& $I_b^* (b\geq 0)$ \\  \hline
$contr_{\nu} (P)$ &$\frac1 2$ & $\frac 3 2$ & $\frac2 3$ & $\frac4 3$&
$\frac{i(b-i)} b$ & $\left \{\begin{array}{ll} 
                 1 , & i=1 \\  
                  1+\frac b 4, & i>1
           \end{array}
\right.$ \\  \hline
$contr_{\nu} (P,Q) (i<j)$& - & - & $\frac1 3$ &$\frac2 3$ &
$\frac{i(b-i)} b$ &  $\left \{\begin{array}{ll} 
                         \frac1 2 , & i=1 \\  
                         \frac{(2+b)}4 ,& i>1
                  \end{array}
\right.  $ \\ \hline
\end{tabular}
\end{center}
}
\end{remark}

\begin{theorem} {\rm (cf.[MP1, Lemma 3.1 and Proposition 3.4])}
 For an elliptic fibration $\pi$: $X \longrightarrow {\bf P}^1$, 
 we have the following formulas:
$$
   \deg J = \sum_{n\geq 1} n(i_n + i^*_n ),
$$ 
and if furthermore $\deg J \not= 0$, then we also have
$$
\deg J \leq  6\sum_{n\geq 1} (i_n + i^*_n ) + 4(ii + iv^* ) +
3(iii + iii^* ) + 2(iv + ii^* ) - 12.
$$ 
\end{theorem}

\section{The possible configurations of   
the unsemi-stable fibrations}

We shall prove Theorem $0.2$ in the present section.

\vskip 0.5cm
\hskip -0.6cm Let $f$: $X \longrightarrow {\bf P}^1$ be a 
(relatively) minimal elliptic surface over ${\bf P}^1$ with
a distinguished section ${\cal O}$. The complete list of 
possible fibers has been given by Kodaira [K1]. It 
encompasses two infinite families $(I_n, I_n^*, n\geq 0)$
and six exceptional cases $(II, III, IV, II^*, III^*, IV^*)$.
And they can be considered as sublattices of the Neron-Severi 
group of $X$ and as such they have rank $(=r(F))$. 
If $e(F)$ denotes the Euler number
of the fiber as a reduced divisor, we can set up the following
table.
\begin{center}
\begin{tabular}{|c|c|c|c|c|c|c|c|c|c|}\hline
 & $I_0$ & $I_n$($n\geq 1$) & $I_n^*$($n\geq 0$)& $II$ &
$III$ & $IV$ & $IV^*$ & $III^*$ & $II^*$ \\  \hline
e & $0$ & $n$ & $n+6$ & $2$ & $3$ & $4$ & $8$ & $9$ & $10$ \\ \hline
r & $0$ & $n-1$ & $n+4$ & $0$ & $1$ &$2$ & $6$ & $7$ & $8$
\\ \hline
\end{tabular} 
\end{center}
\begin{lemma} {\rm (cf.[MP1, Corollary 1.3])}  In all
cases $0\leq e-r \leq 2$. Moreover,

$(1)$ \hskip 0.5cm $e-r=0 \Longleftrightarrow$ the fibre $F$
is smooth, i.e., of type $I_0$;

$(2)$  \hskip 0.5cm $e-r=1 \Longleftrightarrow$ the fibre $F$
is semi-stable, i.e., of type $I_n$ ,$n\geq 1$;

$(3)$ \hskip 0.5cm $e-r=2 \Longleftrightarrow$ the fibre $F$
is unstable.
\end{lemma}
In the following discussion, we denote
\begin{eqnarray*}
{[Q1]} &:=&  \sum_{n\geq 1} ni_n  + \sum_{n\geq 1} (n+6) i^*_n 
+ 6i^*_0 + 10ii^* + 9iii^* + 8iv^* + 4iv + 3iii + 2ii. \\
{[Q2]} &:=& \sum_{n\geq 1} (n-1)i_n  + \sum_{n\geq 1} (n+4) i^*_n 
+ 4i^*_0 + 8ii^* + 7iii^* + 6iv^* + 2iv + iii.  \\
{[Q3]} &:=& \sum_{n\geq 1} i_n  + 2(\sum_{n\geq 1} i^*_n 
+ i^*_0 + ii^* + iii^* + iv^* + iv + iii + ii). \\
{[Q4]} &:=& 6\sum_{n\geq 1} (i_n + i^*_n ) + 4(ii + iv^* ) +
3(iii + iii^* ) + 2(iv + ii^* ) - 12. \\
{[Q5]} &:=&  i^*_0 + iv + iii + ii.
\end{eqnarray*}
It is easy to see $[Q1]=24$, $[Q2]= \rho (X)-2$ and $[Q1]-[Q2]=[Q3]$.
\begin{lemma}
Let $f$: $X \longrightarrow {\bf P}^1$ be an 
elliptic surface over ${\bf P}^1$ with
$\rho (X)=a (\leq 20)$. If $\deg J\not= 0$, then we have

$(1)$ \ \ $[Q4] - \deg J = 6(20 - a - [Q5])$.

$(2)$ \ \ $0 \leq [Q5] \leq 20 - a$. 
\end{lemma}
{\em Proof.}\ By Theorem 1.12, we have
\begin{eqnarray*}
[Q4] - \deg J &=&  [Q4] -  \sum_{n\geq 1} n(i_n + {i}^*_n ) \\
             &=& [Q4] - 24 + \sum_{n\geq 1} 6 i^*_n 
+ 6i^*_0 + 10ii^* + 9iii^* + 8iv^* + 4iv + 3iii + 2ii \\
            &=& 6\sum_{n\geq 1} i_n + 
12( \sum_{n\geq 1}{i}^*_n +  ii^* + iii^* + iv^*) + 6[Q5] - 36 \\
            &=& 6([Q3]-[Q5]-6)\\
            &=& 6(20-a-[Q5]) \geq 0.
\end{eqnarray*}
{\em Proof of Theorem $0.2$.} In this case, $a=20$,
by Lemma 2.2, we get the result.
\begin{lemma} Assume $X$ is an extremal elliptic $K3$ surface 
with an unsemi-stable fibration $f:$ $X \longrightarrow {\bf P}^1$.
Let $m$ be the number of singular fibers of $f$.
Then

$(a)$\  $\sum (e(F) - r(F)) = 6$. 

$(b)$\ $3\leq m \leq 5$.
\end{lemma}
{\em Proof.}\  Since $X$ is $K3$ surface, $\sum e(F)=24$. Also
since $X$ is extremal, $\sum r(F)=\rho -2=18$. This proves (a).
(b) follows from Lemma 2.1 and our definition of {\em unsemi-stable
fibration}.
\vskip 0.5cm
\hskip -0.6cm In the following discussion, we denote $F_i$ (i=1,2,3)
as the singular fiber of $f$ which is not of type $I_n$.
\begin{theorem} 
 Let $f :$ $ X \longrightarrow
{\bf P}^1 $ be an extremal elliptic $K3$ unsemi-stable fibration. Then
the number $m$ of the singular fibers of f is $3$, $4$ or $5$.  

$(A)$ If $m=3$, then the possible fiber types of $(F_1, F_2, F_3)$
are listed in the following.

$(II^* ,I_1^* ,I_1^* )$, $(II^* ,II^* ,IV)$, $(II^* ,IV^* ,I_0^*)$, 
$(III^* ,III^* ,I_0^* )$,

$(III^* ,IV^* ,I_1^* )$, $(III^* ,I_2^* ,I_1^* )$, $(IV^* ,IV^* ,IV^* )$,
$(IV^* ,IV^* ,I_2^* )$,

$(IV^* ,I_2^* ,I_2^* )$, $(IV^* ,I_3^* ,I_1^* )$,$(I_2^* ,I_2^* ,I_2^* )$,
$(I_2^* ,I_3^* ,I_1^* )$, $(I_4^* ,I_1^* , I_1^* )$.

$(B)$ If $m=4$, then the possible fiber types of 
$(F_1, F_2, I_{n_3}, I_{n_4})$ are listed in the following.

$(B.1)$ \ \   $(I_{n_1}^* ,I_{n_2}^* , I_{n_3 } ,I_{n_4 })$ \mbox{where}
$ n_1 \geq n_2 \geq 1$\  \mbox{and}\  $\sum_{i=1}^{4} n_i = 12$.

$(B.2)$          

$(i)$ \ \    $(I_{n_1}^* , II^* ,I_{n_3} ,I_{n_4 })$ 
\mbox{where} $n_1 + n_3 + n_4 =8$\  \mbox{and} \  $n_1 \geq 1$.

$(ii)$ \ $(I_{n_1}^* , III^* ,I_{n_3} ,I_{n_4 })$ 
\mbox{where} $n_1 +n_3 + n_4 =9$\ \mbox{and} \  $n_1 \geq 1$.

$(iii)$\  $(I_{n_1}^* , IV^* ,I_{n_3} ,I_{n_4 })$ 
\mbox{where} $n_1 +n_3 + n_4 =10$ \  \mbox{and}\   $n_1 \geq 1$.

$(B.3.1)$\ 

$(i)$\ $(II^* , II^* ,I_{n_3} ,I_{n_4} )$\  \mbox{ where}\  $n_3 +n_4 =4$.

$(ii)$\ $(II^* , III^* ,I_{n_3} ,I_{n_4} )$\  \mbox{ where}\  $n_3 +n_4 =5$.

$(iii)$\ $(II^* , IV^* ,I_{n_3} ,I_{n_4} )$\  \mbox{ where}\  $n_3 +n_4 =6$.

$(B.3.2)$ \ 

$(i)$\  $(III^* , III^* ,I_{n_3} ,I_{n_4} )$\  \mbox{ where}\  $n_3 +n_4 =6$.

$(ii)$\ $(III^* , IV^* ,I_{n_3} ,I_{n_4} )$\  \mbox{ where}\  $n_3 +n_4 =7$.

$(B.3.3)$\ 

$(IV^* , IV^* ,I_{n_3} ,I_{n_4} )$\  \mbox {where}\  $n_3 +n_4 =8$.

$(C)$ If $m=5$, then the possible fiber types of 
$(F_1, I_{n_2}, I_{n_3}, I_{n_4}, I_{n_5})$ are listed in the following.

$(i)$\ $(I_{n_1}^* ,I_{n_2} ,I_{n_3} ,I_{n_4} ,I_{n_5} )$
\mbox{ where} $\sum_{i=1}^{5} n_i  =18$\  \mbox{and}\  $n_1 \geq 1$.

$(ii)$\ $(II^* ,I_{n_2} ,I_{n_3} ,I_{n_4} ,I_{n_5 } )$\ \mbox{ where}\  
$\sum_{i=2}^{5} n_i =14$.

$(iii)$\ $(III^* ,I_{n_2} ,I_{n_3} ,I_{n_4} ,I_{n_5 } )$\ \mbox{where} \  
$\sum_{i=2}^{5} n_i =15$.

$(iv)$\ $(IV^* ,I_{n_2} ,I_{n_3} ,I_{n_4} ,I_{n_5 } )$\  \mbox{ where} 
\ $\sum_{i=2}^{5} n_i =16$.

\end{theorem}
{\em Proof.}\ We discuss $\deg J = 0$ and  $\deg J\not= 0$
separately.
\vskip 0.5cm
\hskip -0.6cm If $\deg J=0$, then $m=3$. By $[Q2]=18$, we have
$$
18 = iii + 2iv + 6iv^* + 7iii^* +8ii^* + 4i_0^*.
$$
Thus we have the following possible fiber types:
$$
(II^* ,II^* ,IV),\  (II^* ,IV^* ,I_0^*),\  (III^* ,III^* ,I_0^* ),
 \  (IV^* ,IV^* ,IV^* ).
$$
 If $\deg J\not=0$, then by Lemma 2.2, we have
$i^*_0 + iv + iii + ii=0$ and
$$
24 = \sum_{n\geq 1} ni_n  + \sum_{n\geq 1} (n+6) i^*_n +
 10ii^* + 9iii^* + 8iv^* := [a].
$$
If $m=3$, then we have
$$
   24 = \sum_{n\geq 1} (n+6) i^*_n + 10ii^* + 9iii^* + 8iv^* . 
$$
Thus we have the following possible fiber types:
$$
(II^* ,I_1^* ,I_1^* ),\ (III^* ,IV^* ,I_1^* ),
\ (III^* ,I_2^* ,I_1^* ), \ (IV^* ,IV^* ,I_2^* ),\
(IV^* ,I_2^* ,I_2^* ),
$$
$$ 
(IV^* ,I_3^* ,I_1^* ),\   
(I_2^* ,I_2^* ,I_2^* ),\ (I_2^* ,I_3^* ,I_1^* ),\ 
(I_4^* ,I_1^* , I_1^* ).
$$
Combining the above results, we prove  Case(A).
\vskip 0.5cm
\hskip -0.6cm If $m=4$, then with $[Q3] = 6$, we have
$$
  2 = \sum_{n\geq 1}  i^*_n + ii^* + iii^* + iv^* := [b]
$$
and
$$
[a] - 6 \times [b] = \sum_{n\geq 1}ni_n +  \sum_{n\geq 1}ni^*_n + 4ii^* + 3iii^* +
 2iv^* =12. 
$$
This  proves Case(B).
\vskip 0.5cm
\hskip -0.6cm If $m=5$, then with $[Q3] = 6$, we have
$$
 \hskip 2cm   1 = \sum_{n\geq 1}  i^*_n + ii^* + iii^* + iv^* := [c]
$$
and
$$
[a] - 6\times [c] = \sum_{n\geq 1}ni_n +  \sum_{n\geq 1}ni^*_n + 4ii^* + 3iii^* +
 2iv^* = 18. 
$$
This  proves Case(C).

\section {The possible  Mordell-Weil Groups for Case (A) }
 We shall prove the following Theorem 3.1 in the present section. For 
simplity, we label the fiber types which  appeared in 
Case(A) of Theorem 2.4.
\begin{theorem} The possible Mordell-Weil
 Groups for Case $(A)$ are listed in the following table: 
\vskip 0.5cm
\begin{center}
\begin{tabular}{| c|c|c|c|c|c|}\hline
$\sharp$  &  the fibre type &  $MW(f)$ & $\sharp$  & the fibre type & $MW(f)$ \\
 \hline

  $1$ & $(II^* ,I_1^* , I_1^* )$  &  $(0)$   &
 $ 8$  &  $(IV^* ,IV^* , I_2^* )$     &  $(0)$  \\ \hline

 $ 2$ & $ (II^* ,II^* , IV)$ & $(0)$   & 
 $ 9$  &   $(IV^* ,I_3^* , I_1^* )$    &  $(0)$  \\ \hline   

  $3$ & $(II^* ,IV^* , I_0^* )$ & $(0)$    & 
 $ 10$  & $(I_4^* ,I_1^* , I_1^* )$    &  $(0)$, ${\bf Z}/2{\bf Z}$ 
  \\ \hline

$ 4 $& $(III^* ,III^* , I_0^* )$ & $(0)$, ${\bf Z}/2{\bf Z}$  &  
$ 11 $ &  $(I_2^* ,I_2^* , I_2^* )$  &  $(0)$,
${\bf Z}/2{\bf Z}$, 
${\bf Z}/2{\bf Z}\oplus{\bf Z}/2{\bf Z}  $\\ \hline

$ 5$ & $(III^* ,IV^* , I_1^* )$ &  $(0)$     &   
 $12$  &  $(I_2^* ,I_3^* , I_1^* )$     &  $(0)$
  \\ \hline
$ 6$ & $(III^* ,I_2^* , I_1^* )$ & $(0)$, ${\bf Z}/2{\bf Z}$ & 
$13$  & $(IV^* ,I_2^* , I_2^* )$    & $(0)$  \\ \hline
$ 7$ & $(IV^* ,IV^* , IV^* )$   & $(0)$, ${\bf Z}/3{\bf Z}$    &
     &    &   \\ \hline

\end{tabular}  
\end{center}
\end{theorem}
We now explain the outline of the proof of Theorem 3.1. Firstly,
we deal with types 1,2,3,5,8,9,13 (cf. Lemma 3.2). Then we calculate
the possible nontrivial Mordell-Weil Groups of types 4,6,7,10,11 
(cf. Lemma 3.4). Finally we deal with type 12 (cf. Lemma 3.5).
\begin{lemma} For type $m$, 
where $m=1$,$2$,$3$,$5$,$8$,$9$ or $13$, the possible
 Mordell-Weil Group is $(0)$.
\end{lemma}
{\em Proof.}\hskip 0.5cm With Definition 1.9,
Theorem 1.10 and Remark 1.11, it is easy to prove Lemma 3.2. 
For example, $m=1$, if there is a non zero section, say, 
$P_1 \in E(K)_{tor}$, where the $i$-th component of $P_1$
is indicated in Remark 1.11, then by Theorem 1.10, we have  
$$
0=<P_1 ,P_1 > = 2\chi ({\cal O}_X ) + 2(P_1 {\cal O}) 
-\left \{\begin{array}{ll} 
            0 ,    & i=0, \\
          1,     & i=1, \\
          1+\frac{1}{4} , & i>1 .
        \end{array}
\right.
-\left \{\begin{array}{ll}
    0 ,    & i=0, \\
          1,     & i=1, \\
          1+\frac{1}{4} , & i>1 .
        \end{array}
\right.
$$
With $2\chi ({\cal O}_X )=4$ and $(P_1 {\cal O}) \geq 0$, 
we get a contradiction. The others can be proved by the same method.
\begin{remark} {\rm  
In the following calculation,
we let $G_i$, $H_i$ and $J_i$ be the $i$-th
component in the corresponding fiber type $F_1$, $F_2$, $F_3$
respectively (cf.Theorem $2.4$). The numbering of the singular
fiber is defined as following diagrams. Meanwhile, ``$P_1$
pass through the $(i,j,k)$ component'' means $P_1$ only
intersect $G_i$, $H_j$ and $G_k$ in the corresponding fiber
type $F_1$, $F_2$ and $F_3$ respectively. }
\end{remark}
\vskip 0.5cm
\centerline {\font\thinlinefont=cmr5
\begingroup\makeatletter\ifx\SetFigFont\undefined
% extract first six characters in \fmtname
\def\x#1#2#3#4#5#6#7\relax{\def\x{#1#2#3#4#5#6}}%
\expandafter\x\fmtname xxxxxx\relax \def\y{splain}%
\ifx\x\y   % LaTeX or SliTeX?
\gdef\SetFigFont#1#2#3{%
  \ifnum #1<17\tiny\else \ifnum #1<20\small\else
  \ifnum #1<24\normalsize\else \ifnum #1<29\large\else
  \ifnum #1<34\Large\else \ifnum #1<41\LARGE\else
     \huge\fi\fi\fi\fi\fi\fi
  \csname #3\endcsname}%
\else
\gdef\SetFigFont#1#2#3{\begingroup
  \count@#1\relax \ifnum 25<\count@\count@25\fi
  \def\x{\endgroup\@setsize\SetFigFont{#2pt}}%
  \expandafter\x
    \csname \romannumeral\the\count@ pt\expandafter\endcsname
    \csname @\romannumeral\the\count@ pt\endcsname
  \csname #3\endcsname}%
\fi
\fi\endgroup
\mbox
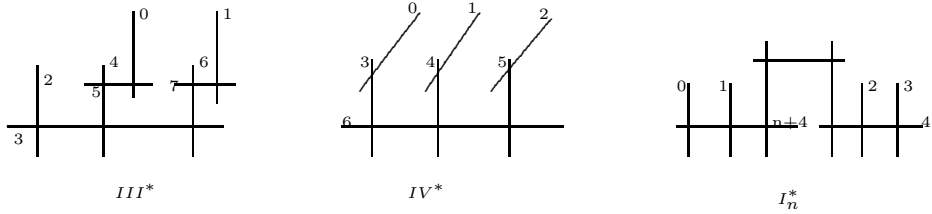
{\beginpicture
\setcoordinatesystem units <0.50000cm,0.50000cm>
\unitlength=0.50000cm
\linethickness=1pt
\setplotsymbol ({\makebox(0,0)[l]{\tencirc\symbol{'160}}})
\setshadesymbol ({\thinlinefont .})
\setlinear
%
% Fig POLYLINE object
%
\linethickness= 0.500pt
\setplotsymbol ({\thinlinefont .})
\plot  1.587 15.875  1.587 15.875 /
%
% Fig POLYLINE object
%
\linethickness= 0.500pt
\setplotsymbol ({\thinlinefont .})
\putrectangle corners at  1.587 15.875 and  7.303 15.875
%
% Fig POLYLINE object
%
\linethickness= 0.500pt
\setplotsymbol ({\thinlinefont .})
\putrectangle corners at  2.381 17.462 and  2.381 15.081
%
% Fig POLYLINE object
%
\linethickness= 0.500pt
\setplotsymbol ({\thinlinefont .})
\putrectangle corners at  4.128 17.462 and  4.128 15.081
%
% Fig POLYLINE object
%
\linethickness= 0.500pt
\setplotsymbol ({\thinlinefont .})
\putrectangle corners at  3.651 16.986 and  5.397 16.986
%
% Fig POLYLINE object
%
\linethickness= 0.500pt
\setplotsymbol ({\thinlinefont .})
\putrectangle corners at  6.032 16.986 and  7.620 16.986
%
% Fig POLYLINE object
%
\linethickness= 0.500pt
\setplotsymbol ({\thinlinefont .})
\putrectangle corners at  6.509 17.462 and  6.509 15.081
%
% Fig POLYLINE object
%
\linethickness= 0.500pt
\setplotsymbol ({\thinlinefont .})
\putrectangle corners at 10.478 15.875 and 16.351 15.875
%
% Fig POLYLINE object
%
\linethickness= 0.500pt
\setplotsymbol ({\thinlinefont .})
\putrectangle corners at 11.271 17.462 and 11.271 15.081
%
% Fig POLYLINE object
%
\linethickness= 0.500pt
\setplotsymbol ({\thinlinefont .})
\putrectangle corners at 14.922 17.621 and 14.922 15.081
%
% Fig POLYLINE object
%
\linethickness= 0.500pt
\setplotsymbol ({\thinlinefont .})
\putrectangle corners at 11.271 17.621 and 11.271 17.462
%
% Fig POLYLINE object
%
\linethickness= 0.500pt
\setplotsymbol ({\thinlinefont .})
\putrectangle corners at  4.921 18.891 and  4.921 16.669
%
% Fig POLYLINE object
%
\linethickness= 0.500pt
\setplotsymbol ({\thinlinefont .})
\putrectangle corners at  7.144 18.891 and  7.144 16.510
%
% Fig POLYLINE object
%
\linethickness= 0.500pt
\setplotsymbol ({\thinlinefont .})
\plot 10.954 16.828 12.541 18.891 /
%
% Fig POLYLINE object
%
\linethickness= 0.500pt
\setplotsymbol ({\thinlinefont .})
\putrectangle corners at 13.018 17.621 and 13.018 15.081
%
% Fig POLYLINE object
%
\linethickness= 0.500pt
\setplotsymbol ({\thinlinefont .})
\plot 12.700 16.828 14.129 18.891 /
%
% Fig POLYLINE object
%
\linethickness= 0.500pt
\setplotsymbol ({\thinlinefont .})
\plot 14.446 16.828 16.034 18.733 /
%
% Fig POLYLINE object
%
\linethickness= 0.500pt
\setplotsymbol ({\thinlinefont .})
\putrectangle corners at 19.367 15.875 and 22.543 15.875
%
% Fig POLYLINE object
%
\linethickness= 0.500pt
\setplotsymbol ({\thinlinefont .})
\putrectangle corners at 23.336 15.875 and 25.876 15.875
%
% Fig POLYLINE object
%
\linethickness= 0.500pt
\setplotsymbol ({\thinlinefont .})
\putrectangle corners at 19.685 16.986 and 19.685 15.081
%
% Fig POLYLINE object
%
\linethickness= 0.500pt
\setplotsymbol ({\thinlinefont .})
\putrectangle corners at 20.796 16.986 and 20.796 15.081
%
% Fig POLYLINE object
%
\linethickness= 0.500pt
\setplotsymbol ({\thinlinefont .})
\putrectangle corners at 21.749 18.098 and 21.749 15.081
%
% Fig POLYLINE object
%
\linethickness= 0.500pt
\setplotsymbol ({\thinlinefont .})
\putrectangle corners at 23.495 18.098 and 23.495 15.081
%
% Fig POLYLINE object
%
\linethickness= 0.500pt
\setplotsymbol ({\thinlinefont .})
\putrectangle corners at 24.289 16.986 and 24.289 15.081
%
% Fig POLYLINE object
%
\linethickness= 0.500pt
\setplotsymbol ({\thinlinefont .})
\putrectangle corners at 25.241 16.986 and 25.241 15.081
%
% Fig POLYLINE object
%
\linethickness= 0.500pt
\setplotsymbol ({\thinlinefont .})
\putrectangle corners at 23.178 15.875 and 23.336 15.875
%
% Fig POLYLINE object
%
\linethickness= 0.500pt
\setplotsymbol ({\thinlinefont .})
\putrectangle corners at 21.431 17.621 and 23.812 17.621
%
% Fig TEXT object
%
\put{\SetFigFont{6}{7.2}{rm}2} [lB] at  2.540 16.986
%
% Fig TEXT object
%
\put{\SetFigFont{6}{7.2}{rm}3} [lB] at  1.746 15.399
%
% Fig TEXT object
%
\put{\SetFigFont{6}{7.2}{rm}4} [lB] at  4.286 17.462
%
% Fig TEXT object
%
\put{\SetFigFont{6}{7.2}{rm}5} [lB] at  3.810 16.669
%
% Fig TEXT object
%
\put{\SetFigFont{6}{7.2}{rm}6} [lB] at  6.668 17.462
%
% Fig TEXT object
%
\put{\SetFigFont{6}{7.2}{rm}7} [lB] at  5.874 16.828
%
% Fig TEXT object
%
\put{\SetFigFont{6}{7.2}{rm}$III^*$} [lB] at  4.445 13.970
%
% Fig TEXT object
%
\put{\SetFigFont{6}{7.2}{rm}0} [lB] at  5.080 18.733
%
% Fig TEXT object
%
\put{\SetFigFont{6}{7.2}{rm}1} [lB] at  7.303 18.733
%
% Fig TEXT object
%
\put{\SetFigFont{6}{7.2}{rm}6} [lB] at 10.478 15.875
%
% Fig TEXT object
%
\put{\SetFigFont{6}{7.2}{rm}0} [lB] at 12.224 18.891
%
% Fig TEXT object
%
\put{\SetFigFont{6}{7.2}{rm}1} [lB] at 13.811 18.891
%
% Fig TEXT object
%
\put{\SetFigFont{6}{7.2}{rm}2} [lB] at 15.716 18.733
%
% Fig TEXT object
%
\put{\SetFigFont{6}{7.2}{rm}3} [lB] at 10.954 17.462
%
% Fig TEXT object
%
\put{\SetFigFont{6}{7.2}{rm}4} [lB] at 12.700 17.462
%
% Fig TEXT object
%
\put{\SetFigFont{6}{7.2}{rm}5} [lB] at 14.605 17.462
%
% Fig TEXT object
%
\put{\SetFigFont{6}{7.2}{rm}0} [lB] at 19.367 16.828
%
% Fig TEXT object
%
\put{\SetFigFont{6}{7.2}{rm}1} [lB] at 20.479 16.828
%
% Fig TEXT object
%
\put{\SetFigFont{6}{7.2}{rm}2} [lB] at 24.448 16.828
%
% Fig TEXT object
%
\put{\SetFigFont{6}{7.2}{rm}3} [lB] at 25.400 16.828
%
% Fig TEXT object
%
\put{\SetFigFont{6}{7.2}{rm}4} [lB] at 25.876 15.875
%
% Fig TEXT object
%
\put{\SetFigFont{6}{7.2}{rm}$IV^*$} [lB] at 12.224 13.970
%
% Fig TEXT object
%
\put{\SetFigFont{6}{7.2}{rm}$I_n^*$} [lB] at 22.066 13.811
%
% Fig TEXT object
%
\put{\SetFigFont{6}{7.2}{rm}n+4} [lB] at 21.907 15.875
\linethickness=0pt
\putrectangle corners at  1.562 19.177 and 25.902 13.811
\endpicture}
} 
\vskip 0.5cm
\begin{lemma} The possible nontrivial Mordell-Weil Group of
type $4$ {\rm (resp. 6,7,10,13)} is ${\bf Z}/2{\bf Z}$
{\rm (resp.  ${\bf Z}/2{\bf Z}$, ${\bf Z}/3{\bf Z}$, 
${\bf Z}/2{\bf Z}$,
${\bf Z}/2{\bf Z}$ or ${\bf Z}/2{\bf Z}\oplus {\bf Z}/2{\bf Z}$)}.
\end{lemma}
{\em Proof.} We only show how to deal with type 6, the others can be
done by the same way. 
\vskip 0.2cm
\hskip -0.6cm  For  m=6, the fiber type is $(III^* ,I_2^* ,I_1^* )$.
 If the  Mordell-Weil Group is nontrivial, then for a non zreo
 section $P_1$, by Theorem 1.10, we have
$$
0=<P_1 ,P_1 > = 2 \chi ({\cal O}_X ) + 2(P_1 O)  
-\left \{\begin{array}{ll}
       0 ,    &i=0,\\
              1 ,    & i=1, \\
          1 + \frac{1}{2} , & i>1, (*) 
        \end{array}
\right.
-\left \{\begin{array}{ll}
     0,             & i=0,\\
            1 ,    & i=1, \\
         1+\frac{1}{2} , & i>1,  (*)
        \end{array}
\right.
-
$$

$$
\left \{\begin{array}{ll}
      0,            & i=0, \\
             1 ,    & i=1, (*)\\
            1+\frac{1}{2} ,    & i>1 .
        \end{array}
\right.
$$
Thus we may assume that the  section $P_1$ pass through
the $(1,2,1)$ component.  An easy  calculation shows
$$
P_1 = {\cal O} + 2F + \sum_{i=1}^{7} \alpha_i G_i + \sum_{j=1}^{6}
\beta_j H_j + \sum_{k=1}^{5} \gamma_k J_k
$$ 
where

\begin{eqnarray*}
  (\alpha_i ) &=& (-\frac3 2 ,-\frac3 2 , -3, -2,-1,-\frac5 2 ,-2), \\ 
 (\beta_j ) &=& (-\frac1 2 ,-\frac3 2 , -1, -2,-\frac3 2 ,-1), \\
 (\gamma_k ) &=& (-1 ,-\frac1 2 ,-\frac1 2 ,-1,-1).
\end{eqnarray*}
and there doesn't exist another non-zero section.
Thus the possible nontrivial Mordell-Weil Group of type $6$ is 
${\bf Z}/2{\bf Z}$.
\begin{lemma} The possible Mordell-Weil Group of type $12$ is
trivial, i.e., $(0)$.
\end{lemma}
{\em Proof.} Assume Lemma 3.5 is false. By the same disscussion
as above,  we may assume that there is a nonzero section $P$ passing 
through the $(1,2,2)$ component.  An easy calculation shows
$$
P = {\cal O} + 2F + \sum_{i=1}^{6} G_i \theta_i + \sum_{j=1}^{7}
\beta_j H_j + \sum_{k=1}^{5} \gamma_k J_k
$$ 
where
\begin{eqnarray*}
  (\alpha_i ) &=& (-1,-\frac1 2 ,-\frac1 2 , -1, -1,-1), \\ 
 (\beta_j ) &=& (-\frac1 2 ,-\frac7 4 , -\frac5 4,-\frac3 2,
 -2,-\frac5 2 ,-1), \\
 (\gamma_k ) &=& (-\frac1 2 ,-\frac5 4 ,-\frac3 4 ,-\frac3 2 ,-1).
\end{eqnarray*}
Thus the possible nontrival Mordell-Weil Group of type $12$ is
${\bf Z}/4{\bf Z}$. That is to say, this group has at least
two nonzero distinct sections, say, $P_1$, $P_2$. On the
other hand, with Theorem 1.10, we have
$$
0=<P_1 ,P_2 > = 2 - (P_1 P_2 )  
- 1
-\left \{\begin{array}{ll}
           1 +\frac{3}{4},    & i=j>1,  \\
         \frac{5}{4} , & j>i>1.  
        \end{array}
\right.
-
\left \{\begin{array}{ll}
             1+\frac{1}{4} ,    & i=j>1, \\
            \frac{3}{4} ,    & j>i>1 .
        \end{array}
\right.
< 0.
$$
Thus we get a contradiction and prove Lemma 3.5.
\vskip 0.5cm
\hskip -0.6cm Combining Lemma 3.2, 3.4 and 3.5,
we prove Theorem 3.1.

\section{The complete determination of the Mordell-Weil Groups for
Case (A) }
 We shall prove Theorem 0.4 in the present section.
\begin{lemma} {\rm (cf.[ATZ, Lemma 3.1])}
Let $S$ be an even symmetric lattice of rank $20$ and 
signature $(1,19)$ and $T$ a positive definite even 
symmetric lattice of rank $2$. Assume that $\varphi$:
$T^{\vee} /T \longrightarrow S^{\vee} /S$ is an
 isomorphism which induces the the following equality 
involving ${\bf Q}/2{\bf Z}$-valued discriminant 
(quadratic) forms:
$$
q_S =-q_T .
$$
Let $X$ be the unique $K3$ surface (up to isomorphisms)
with the transcendental lattice $T_X = T$. Then the Picard
 lattice $Pic X$ is isometric to $S$.
\end{lemma}
\begin{lemma}
 Let $f$: $X \longrightarrow 
{\bf P}^1$ be of type $m$ where  $m=4$,
$6$,$7$,$10$,$11$,$12$ and $13$. Then

$(1)$\   $MW(f_m) \not= (0)$, and further

$(2)$\ $MW(f_{11}) \not= {\bf Z}/2{\bf Z}$.
\end{lemma}
{\em Proof.} Suppose the contrary that  $f$:
 $X \longrightarrow {\bf P}^1$ is of the corresponding type
with $MW(f)=(0)$. Let $(b_{ij} )$ be the intersection metrix
of the transcendental lattice $T=T_X$, then 
$\det{(b_{ij})}=|\det{(PicX)}|$ (cf.[BPV]).  
Modulo congruent action of $SL(2,{\bf Z})$, we may assume that 
$ -b_{11} < 2|b_{12} | \leq b_{11} \leq b_{22} $ , and that 
$b_{12} \geq 0$ when $b_{11} = b_{22} $.

Embed $T$, as a sublattice, into $T^{\vee} = Hom_{\bf Z} (T,
{\bf Z} )$. Then $T^{\vee} /T \cong  (PicX)^{\vee} /(PicX)$.
On the other hand,  $T^{\vee}$ has a {\bf Z}-basis
$(e_1 ,e_2 )(b_{ij} )^{-1} =(g_1 ,g_2 )$, where $e_1$,$e_2$
form a canonical basis of $T$. Then comparing the
 order of ($g_i$) ($i=1$,$2$) with 
 $T^{\vee} /T$, we will get a contraction. For simplity, we only 
show the case for $m=4$, the others can be done by the same way.
\vskip 0.5cm
\begin{center}
\begin{tabular}{|c|c|c|c|}\hline
m & $T^{\vee} /T$ & the possible $T=(b_{ij} )$ & order of $g_1$, $g_2$  \\ \hline
4 & ${\bf Z}/2{\bf Z} \oplus
{\bf Z}/2{\bf Z}  \oplus $
& diag $[2,8]$     & $2$,$8$  \\ \cline{3-4}
  &$( {\bf Z}/2{\bf Z}\oplus {\bf Z}/2{\bf Z})$  & diag $[4,4]$ & $4$,$4$ \\ \hline

\end{tabular}
\end{center}

\begin{remark} {\rm From {\rm Theorem\ 3.1} and
{\rm Lemma\ 4.2}, we know that there does  not exist 
type $12$ or $13$.
The existence of type $m=1$,$3$,$4$,$7$
can be found in {\rm [SI] (page121, 131\ and\ 132)}. Thus the 
Mordell-Weil Group of type $1$ {\rm (resp. $3$,$4$,$7$)}
is $(0)$ {\rm (resp. $(0)$, ${\bf Z}/2{\bf Z}$, 
${\bf Z}/3{\bf Z}$ )}.} 
\end{remark}
\begin{lemma}
Consider the pairs below:
$$
(m, G_m )= (2, (0)),(5, (0)),
(6,{\bf Z}/2{\bf Z} ),(8, (0)), 
 (9, (0)),(10,{\bf Z}/2{\bf Z}), (11,{\bf Z}/2{\bf Z}\oplus {\bf Z}/2{\bf Z}) .
$$
For each of these seven pairs $(m,G_m )$, there is a Jacobian
elliptic $K3$ surface $f_m$: $X_m \longrightarrow {\bf P}^1 $
of type $m$ such that $(m, MW(f_m ))=(m,G_m )$.
\end{lemma}
{\em Proof.} \hskip 0.5cm
 Let $T_m$, $m=2$,$5$,$6$,$8$,$9$, $10$ and $11$ 
be the positive define
symmetric lattice of rank $2$ with the following intersection form,
respectively:
$$
 \left( \begin{array}{rr}
2 & 1 \\
1 & 2
\end{array}
\right) ,  
 \left( \begin{array}{rr}
2 & 0 \\
0 & 12
\end{array}
\right) , 
 \left( \begin{array}{rr}
2 & 0 \\
0 & 4
\end{array}
\right), 
 \left( \begin{array}{rr}
6 & 0 \\
0 & 6
\end{array}
\right) , 
 \left( \begin{array}{rr}
4 & 0 \\
0 & 12
\end{array}
\right) , 
\left( \begin{array}{rr}
4 & 0 \\
0 & 4
\end{array}
\right),
\left( \begin{array}{rr}
2 & 0 \\
0 & 2
\end{array}
\right).
$$
For $m=2$,$5$,$8$,$9$, let $S_m$ be the lattice of rank $20$
and signature $(1,19)$ with the following intersection form, 
respectively
$$
U\oplus  E_8 \oplus E_8 \oplus A_2 ,
U\oplus E_7 \oplus E_6 \oplus D_5 , 
U\oplus E_6 \oplus E_6 \oplus D_6 ,
U\oplus E_6 \oplus D_7 \oplus D_5 .
$$
We now show how to define $S_6$.
 $S_{10}$ and $S_{11}$  can be defined by 
the similar way. 
\vskip 0.5cm
\hskip -0.6cm  Let $\Gamma_6$ be the lattice
$U\oplus E_7 \oplus D_6 \oplus D_5$, with $G_i (1\leq i\leq 7)$, 
$H_j (1\leq j \leq 6)$, $J_k (1\leq k \leq 5)$as the canonical 
basis of $E_7 \oplus D_6 \oplus D_5$ which are indicated in 
{\em Section} 3, and ${\cal O}$, $F$ as a basis of $U$ such
that ${\cal O}^2 =-2$, $F^2 =0$,${\cal O}F=1$.

We extend $\Gamma_6$ to an index-2 integral over lattice
$S_6 = \Gamma_6 + {\bf Z}s_6 $, where
\begin{eqnarray*}
s_6 = {\cal O} + 2F &+& [-\frac3 2 G_1 -\frac3 2 G_2 - 3G_3
-2G_4 -G_5 -\frac5 2 G_6 - 2G_7 ] \\
&+& [-\frac 1 2 H_1 -\frac 3 2 H_2 - H_3 - 2H_4 -\frac 3 2 H_5
     - H_6] \\
&+& [-J_1 - \frac 1 2 J_2 - \frac 1 2 J_3 - J_4 - J_5 ].
\end{eqnarray*}
It is easy to see that intersection form on $\Gamma_6$ can be
extend to an integral even symmetric lattice of signature 
$(1,19)$. Indeed, setting $s=s_6 $, we have
$$
s^2 =-2 , s\cdot F=s\cdot G_1 = s\cdot H_2 = s\cdot J_1 = 1,
$$
$$
 s\cdot G_i = s\cdot H_j =s\cdot J_k =0 \ (\forall i\not=1 ,
j\not= 2 ,k\not= 1 ).
$$

Moreover, $|{\rm det}(S_6 )|= |{\rm det}(\Gamma_6 )| / 2^2  = 8.$
 
Note that $\Gamma_6^{\vee} = {\rm Hom}_{\bf Z} (\Gamma_6 , {\bf Z})$
contains $\Gamma_6$ as a sublettice with $ E_7 \oplus D_6 \oplus D_5$
as the factor group, and is generated by the following, modulo 
$\Gamma_6$:
\begin{eqnarray*}
h_1 &=& (1/2)(G_1 + G_2 + G_6 ), \\
h_2 &=& (1/2)(H_1 + H_2 + H_5 ), \\
h_3 &=& (1/2)(H_1 + H_3 + H_5 ), \\
h_4 &=& (1/4)(2J_1 + J_2 -J_3 + 2J_4 ).
\end{eqnarray*}

Since $(S_6 )^{\vee}$ is an (index-2) sublettice of 
$(\Gamma_6^{\vee})$, an element $x$ is in $(S_6 )^{\vee}$ if
and only if $x = \sum_{i=1}^{4} a_i h_i (\bmod{\Gamma_6} )$ such 
that $x$ is integral on $S_6$, i.e., 
$$
x\cdot s = ( a_1 + a_2 +a_4 )/2
$$
is an integer. Hence $(S_6 )^{\vee}$ is generated by the following
module $\Gamma_6$:
$$
  h_1 + h_2 , h_1 + h_4 ,h_2 +h_4 , h_3 .
$$
Noting that $2h_1$,$2h_2 \in S_6$
and  $ h_1 + h_2  + 2h_4 $ is equal to $s$ (mod $\Gamma_6$)
and hence contained in $S_6$, we find that $(S_6 )^{\vee}$ is
generated by the following, modulo $\Gamma_6$:
$$
 \epsilon_1 = h_3 \ , \ \epsilon_2 =h_1 + h_4. 
$$
Now the fact that $|(S_6 )^{\vee} /S_6 |=8$ and that 
$2\epsilon_1$,$4\epsilon_2 \in S_6 $ imply that
$ (S_6 )^{\vee} /S_6$ is a direct sum of its cyclic
subgroups which are of order $2$,$4$, and generated
by $\epsilon_1$, $\epsilon_2$, modulo $S_6$.

We note that the negative of the discriminant form
$$
-q_{(S_6 )} = (-(\epsilon_1 )^2 )\oplus (-(\epsilon_2 )^2 )
=(3/2) \oplus (3/4).
$$
Similarly, we can get
\begin{eqnarray*}
-q_{(S_{10} )} &=& (-(\epsilon_1 )^2 )\oplus (-(\epsilon_2 )^2 )
=(5/4) \oplus (5/4), \\
-q_{(S_{11} )} &=& (-(\epsilon_1 )^2 )\oplus (-(\epsilon_2 )^2 )
=(1/2) \oplus (1/2),
\end{eqnarray*}
for suitable generators $\epsilon_1$, $\epsilon_2$ in the 
corresponding cases.
\begin{claim}  The pair $(S_m ,T_m )$ satisfies the 
conditions of {\rm Lemma\ 4.1} and hence if we let $X_m$ be the unique
$K3$ surface with $T_{X_m} = T_m$ then ${\rm Pic}X_m = S_m$ 
{\rm (both\ 
two\  equalities\  here\  are\  modulo\  isometries)}.
\end{claim}
{\em Proof of the claim.}\hskip 0.5cm
We need to show that $q_{T_m} = -q_{S_m}$. 

$(1)$\  $m\not =2$. $(S_m )^{\vee} /S_m $ is generated by two
elements $\epsilon_i$ (i=1,2) ($\epsilon_i$ is a simple
sum of the natural generators of $(S_m )^{\vee} /S_m $) such
that for every $a$,$b\in {\bf Z}$ one has $-q_{S_m} (a\epsilon_1
+ b\epsilon_2 ) = -a^2 (\epsilon_1 )^2 - b^2 (\epsilon_2 )^2. $
For all six $m$ where $m\not= 2$, $ \epsilon_i $ can be choosen 
such that $(-\epsilon_1^2 , -\epsilon_2^2 )$ is respectively given
as follows:
\begin{center}
$ (3/2, 7/12)$, $(3/2, 3/4)$, $(5/6, 5/6)$,
$(7/12 , 7/4)$, $(5/4, 5/4)$, $(1/2, 1/2)$.
\end{center}
On the  other hand, $(T_m )^{\vee}$ ($m\not= 2$) 
is generated by
$(g_1, g_2 )=(e_1, e_2 ) T_m^{-1}$, where $e_1$, $e_2$ form a 
canonical basis of $T_m$ which gives rise to the intersection
matrix of $T_m$ shown before this claim. Now the claim follows
from the existence of the following isomorphism,
which induces  $q_{T_m} = -q_{S_m}$:
$$
\phi \ : \ (T_m )^{\vee}/ T_m \longrightarrow   
(S_m )^{\vee}/ S_m , \ \ \ (g_1, g_2 )=(\epsilon_1, \epsilon_2 )B_m.
$$
Here $B_m$ is respectively given as:
$$
  \left( \begin{array}{rr}
1 & 1 \\
6 & 1
\end{array}
\right) ,  
 \left( \begin{array}{rr}
1 & 1 \\
2 & 1
\end{array}
\right) , 
 \left( \begin{array}{rr}
2 & 5 \\
1 & 2
\end{array}
\right) , 
 \left( \begin{array}{rr}
3 & 2 \\
2 & 1
\end{array}
\right) , 
 \left( \begin{array}{rr}
2 & 1 \\
1 & 2
\end{array}
\right),
\left( \begin{array}{rr}
1 & 0 \\
0 & 1
\end{array}
\right)
$$ 
for $m=5$,$6$,$8$,$9$,$10$ and $11$ respectively.
\vskip 0.5cm
\hskip -0.6cm $(2)$ \  $m=2$. In this case, we know $(T_2 )^{\vee}$ 
 is generated by $(g_1, g_2 )=(e_1, e_2 ) T_2^{-1}$,
 where $e_1$, $e_2$ form a 
canonical basis of $T_2$ which gives rise to the intersection
matrix of $T_2$ shown before this claim. In fact we have
$$
  g_1 \equiv g_2 \ ({\rm mod} T_2 ).
$$
Thus $(T_2 )^{\vee} /T_2$ is generated by one element $g_1$,
and the natural isomorphism
$$
   \phi : g_1 \longrightarrow \epsilon_1
$$
will give $q_{T_2} = -q_{S_2}$, where $\epsilon_1 $ is a
canonical {\bf Z}-basis of $A_2$ and 
such that for every $a\in {\bf Z}$, we have
 $-q_{S_2 }(a\epsilon_1 )= -a^2 (\epsilon_1 )^2 $.
\vskip 0.5cm
\hskip -0.6cm  Write $S_m$ (resp. $\Gamma_m$) as 
$U \oplus {\bf B}(m)$ with
${\bf B}(m)$ as in the definitions of them. Let $\cal O$, $F$
be a {\bf Z}-basis of $U$ for all $m$. By [PSS, p.573, Theorem 1],
after an (isometric) action of reflections on $S_m = PicX_m$,
we may assume at the beginning that $F$ is a fibre of elliptic
fibration $f_m$: $X_m \longrightarrow {\bf P}^1$. Since 
${\cal O}^2 = -2$, Riemann-Roch Theorem implies that ${\cal O}$
is an effective divisor for ${\cal O}\cdot F > 0$. Moreover,
${\cal O}\cdot F =1$ implies that ${\cal O}={\cal O}_1 + F'$
where ${\cal O}_1$ is a cross-section of $f_m$ and $F'$ is an 
effective divisor contained in fibres. So $f_m$ is a Jacobian
elliptic fibration and we can choose ${\cal O}_1$ as the zero
element of $MW(f_m )$.

\hskip -0.5cm Let $\Lambda_m $ be the lattice generated by all fiber components
of $f_m$. Clearly, $\Lambda_m = {\bf Z}F \oplus \Delta$, $\Delta
= \Delta (1) \oplus \cdots  \oplus \Delta (r) $ (depending on $m$),
where each $ \Delta (i)$ is a negative definite even
lattice of Dynkin type $A_p$, $D_q$, or $E_r$, contained in a single
 reducible fibre $F_i$ of $f_m$ and spanned by  smooth components
of $F_i$ disjoint from ${\cal O}_1$.

\begin{claim}
We have

$(1)$ \ ${\it Span}_{\bf Z} \{ x\in S_m | x\cdot F=0, x^2 =-2 \} =
\Lambda_m = {\bf Z}F \oplus {\bf B}(m) $; in particular, there
are lattice-isometries: $\Delta \cong {\bf B}(m)$.

$(2)$ \ $MW(f_m) = (0)$ for $m=2$,$5$,$8$ and $9$.

$(3)$  \ $MW(f_m)={\bf Z}/2{\bf Z}$ for $m=6$,$10$.

$(4)$  \ $MW(f_m)={\bf Z}/2{\bf Z}\oplus{\bf Z}/2{\bf Z} $ for $m=11$.
\end{claim}
{\em Proof of the claim.}  The first equality in (1)  
 is from Kodiara's classification of elliptic fibres and the
Riemann-Roch Theorem as used prior to this claim to deduce
${\cal O} \geq 0$. The second equality is clear for that
cases of $m=2$, $5$, $8$ and $9$, because then 
${\rm Pic}X_m = S_m = ({\bf Z}{\cal O} + {\bf Z}F )\oplus
{\bf B}(m)$.

We now show the second equality for $m=6$,$10$ and $11$ using Lemma 1.1. 
Clearly, ${\bf Z}F \oplus {\bf B}(m)$ is contained in the
first term of (1) and hence in $\Lambda_m$. We also have
$$
19={\rm rank}S_m -1 \geq {\rm rank}\Lambda_m  =1+ {\rm rank}\Delta
=  1+ {\rm rank}{\bf B}(m)=19.
$$
Hence $\Delta = \Delta (1) \oplus \cdots \oplus \Delta (r)
\cong \Lambda_m /{\bf Z}F$ contains a finite-index sublettice
$({\bf Z}F\oplus {\bf B}(m))/{\bf Z}F \cong {\bf B}(m)$,

Suppose the contrary that the second equality in (1) is not
true. Then ${\bf B}(m)$ is an index-$n$ $(n>1)$ sublattice 
of $\Delta$.
\vskip 0.5cm
\hskip -0.6cm For $m=6$, by Lemma $1.1$, we know 
$ \Delta =E_7 \oplus D_{11}$. On the other hand, if we
denote $ s_6' = s_6 - {\cal O} - 2F$, then we have 
$$
\Lambda_6  \subset {\it Span}_{\bf Z} \{ x\in S_m | x\cdot F=0 \}=
{\it Span}_{\bf Z} \{ F, G_i , H_j, J_k, s_6' \}
$$
where $1\leq i\leq 7$, $1\leq j\leq 6$ and $1\leq k\leq 5$.

Thus we get ($\bmod{{\bf Z}F}$)
$$
E_7 \oplus D_{11}\subset{\it Span}_{\bf Z} \{ G_i , H_j, J_k, s_6' \}
$$
By a simple  calculation, we find that for any element
$e\in D_{11}-(D_6 \oplus D_5)$, $e\notin {\it Span}_{\bf Z} \{ G_i , H_j, 
J_k, s_6' \}$. Thus we get a contradiction.

Similarly,  we can  prove the second equality in (1) for  
$m=10$, $11$.  The assersion (2),(3) and (4) follow from 
the fact in [Sh, Theorem 1.3],
that $MW(f_m )$ is isomorphic to the factor group of 
Pic$X_m$ modulo $({\bf Z}{\cal O}_1 + {\bf Z}F)\oplus \Delta$,
where the latter is equal to 
$$
({\bf Z}{\cal O} + {\bf Z}F) + \Delta = ({\bf Z}{\cal O} + {\bf Z}F)
\oplus {\bf B}(m) = U \oplus  {\bf B}(m).
$$
This proves the claim. And this completes 
the lattice-theoretical proof of Lemma $4.4$.
\vskip 0.5cm
\hskip -0.6cm Combining Remark 4.3
and Lemma 4.4, we prove Theorem 0.4.

\section*{References}

\hskip 0.6cm [ATZ]\ E.Artal Bartolo, H.Tokunaga and D.Q.Zhang:
{\em Miranda-Persson's Problem on Extremal Elliptic $K3$ Surfaces},
{\bf alg-geom/9809065}.

[BPV]\ W.P.Barth, C.A.M.Peters and A.J.H.M.Van de Ven:
{\em Compact complex surfaces}, Springer, Berlin, 1984.

[CS]\ J.Conway and N.Sloane:{\em  Sphere Packings, Lattices and
Groups}, Grund.Math.Wiss.{\bf 290}, Springer-Verlag (1988). 

[K1]\ K.Kodaira: {\em  On compact complex 
analytic surfaces II}, Ann.Math. {\bf 77}(1963), 563-626.

[MP1]\ R.Miranda and U.Persson: {\em On extremal rational elliptic 
surfaces}, Math.Z.{\bf 193}(1986), 537-558.

[MP2]\ R.Miranda and U.Persson: {\em Configurations of $I_n$ fibers
on elliptic $K3$ surfaces}, Math.Z.{\bf 201}(1989), 339-361.

[MP3]\ R.Miranda and U.Persson: {\em Mordell-Weil Groups of extremal
elliptic $K3$ surfaces}, Problems in the theory of
surfaces and their classification (Cortona, 1988), Symposia 
Mathematica, XXXII, Academic Press, London,1991, 167-192.

[PSS]\ I.I.Pjateckii-Sapiro and I.R.Safarevic: {\em Torelli's 
theorem for algebraic surfaces of type $K3$}, Math.USSR Izv.
{\bf 5} (1971), 547-588.

[SZ]\ I.Shimada and D.Q.Zhang : {\em Classification of extremal
elliptic $K3$ surfaces and fundamental groups of open $K3$ surfaces},
preprint,1999.

[Sh]\ T.Shioda: {\em  On the Mordell-Weil Lattices}, Commentarii
Mathematici, Universitatis Sancti Pauli, Vol.{\bf 39}(1990),
211-240.

[SI]\ T.Shioda and H.Inose: {\em On singular $K3$ surfaces},
 Complex analysis and algebraic geometry, papers dedicated to K.Kodiaira,
Iwanami Shoten and Cambridge University Press, London, 1977,
119-136.

\end{document}